\documentclass{article}

\usepackage{amsfonts}
\usepackage{amssymb}
\usepackage{epsfig}

\newcommand{\N}{{\mathbb N}}
\newcommand{\Z}{{\mathbb Z}}
\newcommand{\X}{{\bar X}}
\newcommand{\Y}{{\bar Y}}
\newcommand{\U}{{\bar U}}
\newcommand{\V}{{\bar V}}

\newcommand{\g}{{\bar g}}
\newcommand{\h}{{\bar h}}

\begin{document}
\begin{center}
\large
ON PRODUCTS OF QUASICONVEX SUBGROUPS IN HYPERBOLIC GROUPS

\vspace{.5cm}
Ashot Minasyan
\end{center}


\begin{abstract}
An interesting question about quasiconvexity in a hyperbolic group concerns finding
classes of quasiconvex subsets that are closed under finite intersections.
A known example is the class of all quasiconvex subgroups [1].
However, not much is yet learned about the structure of arbitrary quasiconvex subsets.
In this work we study the properties of products of quasiconvex subgroups; we show
that such sets  are quasiconvex, their finite intersections have a similar algebraic representation
and, thus, are quasiconvex too.
\end{abstract}

\vspace{.5cm}
{\noindent \bf 0. Introduction}

\vspace{.15cm}
Let $G$ be a hyperbolic group, $\Gamma(G,{\cal A})$ -- its Cayley graph corresponding
to a finite symmetrized generating set $\cal A$ (i.e. for each element $a \in \cal A$, $a^{-1}$ also belongs to this set).
 A subset $Q \subseteq G$ is said to be $\varepsilon$-{\it quasiconvex},
 if any geodesic connecting two elements from $Q$
belongs to a closed $\varepsilon$-neighborhood ${\cal O}_{\varepsilon}(Q)$ of $Q$ in $\Gamma(G,{\cal A})$
for some $\varepsilon \ge 0$. $Q$ will be called {\it quasiconvex} if there exists $\varepsilon > 0$ for which
it is $\varepsilon$-quasiconvex.

In [4] Gromov proves that the notion of quasiconvexity in a hyperbolic group
does not depend on the choice of a finite generating set (it is easy to show that this is not true in
an arbitrary group).

If $A,B \subseteq G$ then their product is a subset of $G$ defined by \\ $A\cdot B = \{ ab~|~a \in A,~b\in B\} $.

\vspace{.15cm}
{\bf \underline{Proposition 1.}} If the sets $A_1, \dots,A_n \subset G$ are quasiconvex
then their pro\-duct set $A_1A_2 \cdot \dots\cdot A_n \stackrel{def}{=} \{
a_1a_2 \cdot \dots \cdot a_n~|~a_i \in G_i \} \subset G$ is also quasiconvex.

\vspace{.15cm}
Proposition 1 was proved by Zeph Grunschlag in 1999 in [11; Prop. 3.14] and, independently,
by the author in his diploma paper in 2000.

If $H$ is a subgroup of $G$ and $x \in G$ then the subgroup conjugated to $H$ by $x$ will be denoted
$H^x = xHx^{-1}$. The main result of the paper is

\vspace{.15cm}
{\bf \underline{Theorem 1.}}
Suppose $G_1, \dots, G_n$, $H_1, \dots, H_m$ are quasiconvex sub\-groups
of the group $G$, $n,m \in \N$; $f,e \in G$.
Then there exist numbers \\ $r, t_l \in \N \cup \{0\}$ and $f_l,\alpha_{lk},
\beta_{lk} \in G$, $k =1,2,\dots, t_l$ (for every
fixed $l$), \\$l =1,2,\dots, r$, such that
$$fG_1 G_2 \cdot \dots \cdot G_n \cap eH_1 H_2
\cdot \dots \cdot H_m = \bigcup_{l=1}^r f_l S_l $$
where for each $l$, $t = t_l$, there are indices $1 \le i_1\le i_2\le \dots \le i_t \le n$,
$1 \le j_1\le $ $ \le j_2\le \dots \le j_t \le m$~:
$$ \quad S_l = (G_{i_1}^{\alpha_{l1}} \cap H_{j_1}^{\beta_{l1}}) \cdot
\dots \cdot (G_{i_t}^{\alpha_{lt}} \cap H_{j_t}^{\beta_{lt}}). $$

\vspace{.15cm}
This claim does not hold if the group $G$ is not hyperbolic~: set
$G_1 = \langle (1,0) \rangle$, $G_2 = \langle (0,1) \rangle$,
$H = \langle (1,1) \rangle $ -- cyclic subgroups of $ \Z^2$ (they are quasiconvex in
$\Z^2$ with generators $\{(1,0), (0,1),(1,1) \}$)~. $G_1\cdot G_2 = \Z^2$, thus,
$G_1G_2 \cap H = H$ but $G_1 \cap H = G_2\cap H = \{(0,0)\}$ and, if the
statement of the theorem held for $\Z^2$ then $H$ would be finite -- a contradiction.

The above example can also be used as another argument to prove the well-known fact
that $\Z^2$ can not be embedded into a hyperbolic group (because any cyclic subgroup is
quasiconvex in a hyperbolic group).

The condition that the subgroups $G_i$, $H_j$ are quasiconvex is also necessary: using Rips'
Construction ([12]) one can achieve a group $G$ satisfying the small cancellation condition
$C'(1/6)$ (and, therefore, hyperbolic) and its finitely generated normal subgroup $K$ such that
$G/K \cong \Z^2$. Let $\phi$ be the natural epimorphism from $G$ to $\Z^2$,
$G_1 = \phi^{-1}(\langle (1,0)\rangle) \le G$, $G_2 = \phi^{-1}(\langle (0,1)\rangle) \le G$,
$H = \phi^{-1}(\langle (1,1)\rangle) \le G$. $G_1,G_2$ and $H$ are finitely generated subgroups of
$G$, $G_1 \cdot G_2 = G$ because $\langle (1,0)\rangle \cdot \langle (0,1)\rangle = \Z^2$
and $K \le G_2$, thus $G_1 \cdot G_2 \cap H = H$. But for every $\alpha, \beta \in G$
$\phi (G_i^\alpha \cap H^\beta) \subseteq \phi(G_i)^{\phi(\alpha)} \cap \phi (H)^{\phi(\beta)} = \{(0,0)\}$,
$i=1,2$. Hence, it is impossible to obtain the infinite subgroup $\phi (H)$ from products of cosets to such sets, and
we constructed the counterexample needed.

\vspace{.15cm}
{\bf Definition :} let $G_1,G_2,\dots,G_n$ be quasiconvex subgroups of $G$,
$f_1,f_2,\dots,$ $f_n \in G$, $n \in \N$. Then the set
$$f_1G_1f_2G_2 \cdot \dots\cdot f_nG_n = \{ f_1g_1f_2g_2 \cdot \dots \cdot f_ng_n \in G~|~g_i \in G_i,~i=1,\dots,n\}
$$ will be called {\it quasiconvex product}.

\vspace{.15cm}
{\bf \underline{Corollary 2.}} An intersection of finitely many quasiconvex products
is a finite union of quasiconvex products.

\vspace{.15cm}
Thus the class of finite unions of quasiconvex products is closed under taking finite
intersections.

Recall that a group $H$ is called {\it elementary} if it has a
cyclic subgroup $\langle h \rangle$ of finite index. An elementary subgroup
of a hyperbolic group is quasiconvex (see remark 5, Section 4).
It is well known that any element $x$ of infinite order in $G$ is contained in a
unique maximal elementary subgroup $E(x) \leqslant G$ [4], [5].  \\
Every non-elementary hyperbolic group contains the free group of rank 2 [5, Cor. 6].

Suppose $G_1,G_2,\dots,G_n,H_1,H_2,\dots,H_m$ are infinite maximal elementary \\
subgroups of $G$, $f,e \in G$. And $G_i \neq G_{i+1}$, $H_j \neq H_{j+1}$, $i=1,\dots,n-1$, $j=1,\dots,m-1$.
Then we present the following uniqueness result for the products of such subgroups:

\vspace{.15cm}
{\bf \underline{Theorem 2.}} The sets $fG_1\cdot \dots \cdot G_n$ and $eH_1\cdot \dots \cdot H_m$
are equal if and only if $n=m$, $G_n=H_n$, and there exist elements $z_j \in H_j$, $j=1,\dots,n$, such that
$G_j = (z_{n} z_{n-1} \dots z_{j+1}) \cdot H_j \cdot
(z_{n} z_{n-1} \dots z_{j+1})^{-1}$, $j=1,2,\dots,n-1$, $f = ez_1^{-1}z_2^{-1}\dots z_n^{-1}$.

\vspace{.15cm}
Similarly to quasiconvex products one can define ME-products to be the products
of cosets of maximal elementary subgroups in $G$ (the full definition is given in Section 4).
The statement of Corollary 2 can be strengthened in this case :

\vspace{.15cm}
\underline{\bf Theorem 3.} Intersection of any family (finite or infinite) of finite unions of ME-products is
 a finite union of  ME-products.

\vspace{.15cm}
An example which shows that an analogous property is not true for arbitrary quasiconvex products is constructed at the
end of this paper.

Thus, all finite unions of ME-products constitute a topology ${\cal T}$ (of closed sets) on the set of elements
of a hyperbolic group. Taking an inverse, left and right shifts in $G$ are continuous operations in ${\cal T}$.
Also, by definition, any point is closed in ${\cal T}$, so ${\cal T}$ is weakly separated ($T_1$).
However, if $G$ is infinite elementary, then ${\cal T}$ turns out to be the topology of finite complements which is
not Hausdorff, also, in this case, the group multiplication is not continuous with respect to  ${\cal T}$ (since
any product of two  non-empty open sets contains the identity of $G$).

\vspace{.5cm}
{\noindent \bf 1. Preliminary information}

\vspace{.15cm}
Let $d(\cdot,\cdot)$ be the usual left-invariant metric on the Cayley graph of the group $G$ with generating set
$\cal A$. For any two points $x,y \in \Gamma(G,{\cal A})$ we fix a geodesic path between them and denote it by $[x,y]$.

If $Q \subset \Gamma(G,{\cal A})$, $N \ge 0$, the closed $N$-neighborhood of $Q$ will be denoted by
$${\cal O}_N (Q) \stackrel{def}{=} \{x\in \Gamma(G,{\cal A})~|~d(x,Q) \le N \}~.$$

If $x,y,w \in \Gamma(G,{\cal A})$, then the number
$$(x|y)_w \stackrel{def}{=} \frac12 \Bigl(d(x,w)+d(y,w)-d(x,y) \Bigr)$$
is called the {\it Gromov product} of $x$ and $y$ with respect to $w$.

Let $abc$ be a geodesic triangle. There exist "special" points $O_a \in [b,c]$,\\ $O_b \in [a,c]$, $O_c \in [a,b]$ with the properties:
$d(a,O_b) = d(a,O_c) = \alpha$, $d(b,O_a) = $ $=d(b,O_c) = \beta$, $d(c,O_a) = d(c,O_b) = \gamma$. From a corresponding
system of linear equations one can find that $\alpha = (b|c)_a$, $\beta = (a|c)_b$, $\gamma = (a|b)_c$. Two points \\
$O \in [a,b]$ and $O' \in [a,c]$ are called $a$-{\it equidistant} if $d(a,O) = d(a,O') \le \alpha$.  \\
The triangle $abc$ is said to be $\delta$-{\it thin} if for any two points $O,O'$ lying on its sides and
equidistant from one of its vertices, $d(O,O') \le \delta$ holds. \\
$abc$ is $\delta$-{\it slim} if each of its sides belongs to a closed $\delta$-neighborhood of the two others.

We assume the following equivalent definitions of hyperbolicity of $\Gamma(G,{\cal A})$
to be known to the reader (see [6], [2]): \\
$1^{\circ}.$ There exists $\delta \ge 0$ such that for any four points $x,y,z,w \in \Gamma(G,{\cal A})$ their Gromov products satisfy
$$(x|y)_w \ge min\{(x|z)_w,(y|z)_w\} - \delta~;$$
$2^{\circ}.$ All triangles in $\Gamma(G,{\cal A})$ are {\cal $\delta$-thin} for some $\delta \ge 0$;\\
$3^{\circ}.$ All triangles in $\Gamma(G,{\cal A})$ are $\delta$-slim for some $\delta \ge 0$.

Now and below we suppose that $G$ meets $1^{\circ},2^{\circ}$ and $3^{\circ}$ for a fixed (sufficiently large)
$\delta \ge 0$. $3^{\circ}$ easily implies

\vspace{.15cm}
\underline{\bf Remark 0.} Any side of a geodesic $n$-gon ($n \ge 3$) in $\Gamma(G,{\cal A})$
belongs to a closed $(n-2)\delta$-neighborhood of the union of the rest of its sides.

\vspace{.15cm}
Let $p$ be a path in the Cayley graph of $G$. Further on by $p_-$, $p_{+}$ we will denote the startpoint and
the endpoint of $p$, by $||p||$ -- its length; $lab(p)$, as usual, will mean the word in the alphabet $\cal A$
written on $p$. $elem(p) \in G$ will denote the element
of the group $G$ represented by the word $lab(p)$.

A path $q$ is called $(\lambda,c)$-quasigeodesic if there exist $0<\lambda \le 1$, $c \ge 0$, such that
for any subpath $p$ of $q$ the inequality $\lambda ||p|| - c \le d(p_-,p_+)$ holds.\\
In a hyperbolic space quasigeodesics and geodesics with same ends are mutually close~:

\vspace{.15cm}
\underline{\bf Lemma 1.1.} ([6; 5.6,5.11], [2; 3.3]) There is a constant $N=N(\delta,\lambda,c)$ such that for any
$(\lambda,c)$-quasigeodesic path $p$ in $\Gamma(G,{\cal A})$ and a geodesic
$q$ with $p_- = q_-$, $p_+ = q_+$, one has $p \subset {\cal O}_N(q)$ and $q \subset {\cal O}_N(p)$.

\vspace{.15cm}
An important property of cyclic subgroups in a hyperbolic group states

\vspace{.15cm}
\underline{\bf Lemma 1.2.} ([6; 8.21], [2; 3.2]) For any word $w$ representing an element $g \in G$ of infinite order  there
exist constants $\lambda >0$, $c \ge 0$, such that any path with a label $w^m$ in the Cayley graph of
$G$ is $(\lambda,c)$-quasigeodesic for arbitrary integer $m$.

\vspace{.15cm}
A broken line $p=[X_0,X_1,\dots,X_n]$ is a path obtained as a consequent concatenation of geodesic segments
$[X_{i-1},X_{i}]$, $i=1,2,\dots,n$.
Later, in this paper, we will use the following fact concerning broken lines in a hyperbolic space:

\vspace{.15cm}
\underline{\bf Lemma 1.3.} ([3, Lemma 21]) Let $p = [X_0,X_1,\dots,X_n]$ be a broken line in $\Gamma(G,{\cal A})$ such that
$||[X_{i-1},X_i]|| > C_1$ $\forall~i=1,\dots,n$, and $(X_{i-1}|X_{i+1})_{X_i} \le C_0$ $\forall~i=1,\dots,n-1$,
where $C_0 \ge 14 \delta$, $C_1 > 12(C_0 + \delta)$. Then $p$ is contained in the closed $2C_0$-neighborhood
${\cal O}_{2C_0}([X_0,X_n])$ of the geodesic segment $[X_0,X_n]$.

\vspace{.15cm}
Suppose $H= \langle {\cal X} \rangle$ is a subgroup of $G$ with a finite symmetrized generating set $\cal X$.
If $h \in H$, then by $|h|_G$ and $|h|_H$ we will denote the lengths of the element $h$ in $\cal A$ and
$\cal X$ respectively.
The {\it distortion function} $D_H: ~\N \to \N$ of $H$ in $G$ is defined by $D_H(n) = max\{|h|_H~|~h \in H, ~|h|_G \le n \}$. \\
If $\alpha, \beta~: \N \to \N$ are two functions then we write $\alpha \preceq \beta$ if $\exists~K_1,K_2>0$~:
$\alpha(n) \le K_1\beta(K_2n)$. $\alpha(n)$ and $\beta(n)$ are said to be equivalent if $\alpha \preceq \beta$ and
$\beta \preceq \alpha$.

Evidently, the function $D_H$ does not depend (up to this equivalence) on the choice of finite generating sets $\cal A$ of $G$
and $\cal X$ of $H$. One can also notice that $D_H(n)$ is always at least linear (provided that $H$ is infinite).
If $D_H$ is equivalent to linear, $H$ is called {\it undistorted}.

\vspace{.15cm}
\underline{\bf Lemma 1.4.} ([2; 3.8], [7; 10.4.2]) A quasiconvex subgroup $H$ of a hyperbolic group $G$ is finitely
generated.

\vspace{.15cm}
{\bf \underline{Remark 1.}} From the proof of this statement it also follows that $D_H$ is equivalent to linear for a quasiconvex subgroup $H$.

\vspace{.15cm}
Indeed, it was observed in [2] that
if $H$ is $\varepsilon$-quasiconvex, it is generated by
finitely many elements $x_i,~i = 1, \dots , s$, such that $|x_i|_G \le
2\varepsilon+1$ $\forall~ i$, and $\forall~ h \in H$, $h =a_1\cdot \dots \cdot a_r$,
$a_j \in {\cal A}$, hence
$\exists~ i_1, \dots,i_r \in \{1,2,\dots,s\}$: $h = x_{i_1}x_{i_2}\cdot \dots \cdot x_{i_r}$.

The proof of corollary 2 is based on

\vspace{.15cm}
\underline{\bf Lemma 1.5.} ([1; Prop. 3]) Let $G$ be a group generated by a finite set $\cal A$. Let $A,B$ be
subgroups of $G$ quasiconvex with respect to $\cal A$. Then $A \cap B$ is quasiconvex with respect to $\cal A$.

\vspace{.15cm}
We will use the following notion in this paper~:

\vspace{.15cm}
{\bf Definition :} let $H = \langle {\cal X} \rangle \le
G = \langle {\cal A} \rangle$, $card({\cal X}) < \infty$, $card({\cal A}) < \infty$.
A path $P$ in $\Gamma(G,{\cal A})$ will be called $H$-{\it geodesic}
(or just $H$-{\it path}) if~: \\
a) $P$ is labelled by the word  $ a_{11} \dots a_{1k_1} \dots a_{s1}
\dots a_{sk_s} $  corresponding to an element $elem(P)=x \in H$, where $a_{ij} \in {\cal A}$; \\
b) $a_{j1} \dots a_{jk_j}$  is a shortest word for generator $x_j \in {\cal X}$ (i.e. $|x_j|_G = k_j$),
$j=1, \dots, s$~; \\
c) $x=x_1 \dots x_s$ in $H$, $|x|_H = s$. \\
I.e. $P$ is a broken line in $\Gamma(G,{\cal A})$ with segments
corresponding to shortest representations of generators of $H$ by means
of $\cal A$.

\vspace{.15cm}
{\bf \underline{Lemma 1.6.}} (see also [10; Lemma 2.4]) Let $H$ be a (finitely generated) subgroup of a
$\delta$-hyperbolic group $G$. Then $H$ is quasiconvex iff $H$ is undistorted in $G$.

\vspace{.15cm}
{\it Proof.} The necessity is given by remark 1.

To prove the sufficiency, suppose $H=\langle {\cal X} \rangle $,
$card({\cal X}) < \infty$, and
$D_H(n) \le cn$, $\forall~n \in \N$, $c>0$. For arbitrary two vertices $x,y \in H$ there is a $H$-path $q$ connecting them in
$\Gamma(G,{\cal A})$. Let $p$ be any its subpath. By definition, there exists a subpath $p'$ of $q$ such
that $p'_{-},p'_{+} \in H$, subpaths of $q$ from $p_{-}$ to $p'_{-}$ and from $p_{+}$ to $p'_{+}$
are geodesic, and $d(p_{-},p'_{-}) \le \varkappa/2$, $d(p_{+},p'_{+}) \le \varkappa/2$, where
$\varkappa = max\{|h|_G~|~h \in {\cal X}\} < \infty$. In particular, $p'$ is also $H$-geodesic.\\
Using the property c) from the definition of a $H$-path we obtain
$$||p'|| \le \varkappa \cdot |elem(p')|_H \le \varkappa \cdot c \cdot d(p'_{-},p'_{+})~.$$
Therefore, $||p|| \le ||p'|| +\varkappa \le \varkappa \cdot c \cdot d(p'_{-},p'_{+})  + \varkappa \le
\varkappa \cdot c \cdot d(p_{-},p_{+})  + \varkappa^2 c +\varkappa$, which shows that $q$ is
$(\frac1{\varkappa c},\varkappa+\frac1c)$-quasigeodesic. By lemma 1.1 $\exists~N=N(\varkappa,c) $ such that
any geodesic path between $x$ and $y$ belongs to the closed $N$-neighborhood ${\cal O}_N(q)$
but $q \subset {\cal O}_{\varkappa/2}(H)$ in the Cayley graph of $G$. Hence, $H$ is
quasiconvex with the constant $(N+\varkappa/2)$, and the lemma is proved. $\square$

\vspace{.15cm}
During this proof we showed

\vspace{.15cm}
{\bf \underline{Remark 2.}} If $H$ is a quasiconvex subgroup of a hyperbolic
group $G$ then any $H$-path is $(\lambda,c)$-quasigeodesic for some
$\lambda,c$ depending only on the subgroup $H$.



\vspace{.15cm}
Let the words $W_1,\dots,W_l$ represent elements $w_1,\dots,w_l$ of infinite order in a hyperbolic group $G$.
For a fixed constant $K$  consider the set \\ $S_M = S(W_1,\dots,W_l;K,M)$ of words
$$W = X_0W_1^{\alpha_1}X_1W_2^{\alpha_2}X_2 \dots W_l^{\alpha_l}X_l$$
where $||X_i|| \le K$ for $i=0,1,\dots,l$, $|\alpha_2|,\dots,|\alpha_{l-1}| \ge M$, and the element of $G$ represented by
$X_i^{-1}W_iX_i$ does not belong to the maximal elementary subgroup $E(w_{i+1}) \le G$ containing $w_{i+1}$ for $i=1,\dots,l-1$.

\vspace{.15cm}
\underline{\bf Lemma 1.7.} ([5; Lemma 2.4]) There exist constants $\lambda>0$, $c \ge 0$ and \\ $M>0$ (depending on $K$,
$W_1,\dots,W_l$) such that any path in $\Gamma(G,{\cal A})$ labelled by an arbitrary word $W\in S_M$ is
$(\lambda,c)$-quasigeodesic.
\vspace{.15cm}

\vspace{.15cm}
\underline{\bf Lemma 1.8.} Suppose $l \in \N$, $K>0$, and $w_1,\dots,w_l \in G$
are elements of infinite order. Then there are $\lambda>0$, $c \ge 0$ and $M>0$ (depending on $K$,
$w_1,\dots,w_l$) such that for arbitrary $x_0,x_1,\dots,x_l \in G$, $|x_i|_G \le K$, $i=0,\dots,l$,
with conditions~~  $w_i \notin x_iE(w_{i+1})x_i^{-1}$ $\forall~i \in \{1,\dots,l-1\}$, and any
$\alpha_i \in \Z$, $|\alpha_i| \ge M$, $i=2,\dots,l-1$, the element
$$w = x_0w_1^{\alpha_1}x_1w_2^{\alpha_2}x_2 \cdot \dots \cdot w_l^{\alpha_l}x_l \in G $$
satisfies $|w|_G \ge \lambda |\alpha_1| -c$.

\vspace{.15cm}
{\it Proof.} As follows from Lemma 1.7 and the definition of a $(\lambda,c)$-quasigeodesic path, one has the following inequality:
$$|w|_G \ge \lambda \cdot \left( |x_0|_G + \sum_{i=1}^l (|\alpha_i| |w_i|_G+|x_i|_G) \right) -c \ge \lambda \cdot |\alpha_1| |w_1|_G -c
\ge \lambda |\alpha_1| -c~. \quad \square$$

\vspace{.5cm}
{\noindent \bf 2. Quasiconvex sets and their products}

\vspace{.15cm}
{\bf \underline{Remark 3.}} Any finite subset of $G$ is $d$-quasiconvex
(where $d$ is the diameter of this set).

\vspace{.15cm}
{\bf \underline{Remark 4.}}  Let $Q\subseteq G$ be $\varepsilon$-quasiconvex, $g\in G$. Then
(a)~the left shift $gQ=\{ gx~|~x\in Q\}$ is quasiconvex with the same constant;
(b)~the right shift $Qg=\{ xg~|~x\in Q\}$ is quasiconvex (possibly, with a different quasiconvexity constant).

\vspace{.15cm}
(a) holds because the metric $d(\cdot,\cdot)$ is left-invariant. \\
$x,y \in Q$ if and only if $xg,yg \in Qg$. By remark 0
$$[xg,yg] \subset {\cal O}_{2\delta} \Bigl( [x,xg] \cup [x,y] \cup [y,yg] \Bigr) \subset
{\cal O}_{2\delta +|g|_G} \Bigl( [x,y]\Bigr) \subset {\cal O}_{2\delta +|g|_G+\varepsilon}(Q) \subset$$
$$\subset  {\cal O}_{2\delta +2|g|_G+\varepsilon}(Qg)$$ therefore (b) is true.

\vspace{.15cm}
Therefore, a left coset of a quasiconvex subgroup and a conjugate subgroup to it are quasiconvex
(in a hyperbolic group).

\vspace{.15cm}
{\bf \underline{Lemma 2.1.}} (see also [11; Prop. 3.14]) A finite union of quasiconvex sets in a hyperbolic group
$G$ is quasiconvex.

\vspace{.15cm}
{\it Proof.} It is enough to prove that if $A,B \subset G$ are $\varepsilon_i$-quasiconvex, $i=1,2$, respectively,
then $C=A\cup B$ is quasiconvex. \\
If both $x,y \in A$ or $x,y \in B$ then $[x,y] \subset {\cal O}_{max\{\varepsilon_1,\varepsilon_2\}}(C)$. So, assume that
$x \in A$, $y \in B$. Fix $a \in A$, $b \in B$, and consider the geodesic quadrangle $xyba$ (see Figure 1).

\begin{figure}
\begin{center}
\input{fig1.tex}

Figure 1
\end{center}

\end{figure}

By remark 0 we have $[x,y] \subset {\cal O}_{2\delta}([x,a] \cup [a,b] \cup [b,y])$. After denoting $d(a,b) = 2\eta$ we obtain
$[x,a] \subset {\cal O}_{\varepsilon_1}(A)$, $[b,y] \subset {\cal O}_{\varepsilon_2}(B)$,
$[a,b] \subset {\cal O}_{\eta}(A\cup B)$. Hence $[x,a] \cup [a,b] \cup [b,y] \subset {\cal O}_{max\{\varepsilon_1,\varepsilon_2,\eta \}}(C)$,
consequently, \\ $[x,y] \subset {\cal O}_{max\{\varepsilon_1,\varepsilon_2,\eta\}+2\delta}(C)$, and the lemma is proved.
$\square$


\vspace{.15cm}
{\it Proof} of Proposition 1. Assume $n=2$ (for $n>2$ the statement will follow by induction). \\
So, let $A$, $B$ be $\varepsilon_i$-quasiconvex subsets of $G$ respectively, $i=1,2$.

Consider arbitrary $a_1b_1,a_2b_2 \in AB$, $a_i \in A$, $b_i \in B_i$, $i=1,2$, and fix an
element $b \in B$, $|b|_G = \eta$. Then, since the triangles are $\delta$-slim,
$$[b_1,1_G] \subset {\cal O}_{\delta} ([b,1_G] \cup [b,b_1]) \subset
{\cal O}_{\delta+\eta}([b,b_1]) \subset {\cal O}_{\delta+\eta+\varepsilon_2}(B)~.$$
Denoting $\varepsilon_3 = \delta+\eta+\varepsilon_2$, one obtains
$[b_1,1_G] \subset {\cal O}_{\varepsilon_3} (B)$ and, similarly, \\ $[b_2,1_G] \subset {\cal O}_{\varepsilon_3} (B)$.
Therefore, $[a_1b_1,a_1] \subset {\cal O}_{\varepsilon_3} (a_1B)$, $[a_2b_2,a_2] \subset {\cal O}_{\varepsilon_3} (a_2B)$.
Also, observe that $\forall ~a \in A$ ~~$d(a,ab) = |b|_G = \eta$, i.e. $A\subset {\cal O}_{\eta}(Ab)
\subset {\cal O}_{\eta}(AB)$, hence $[a_1,a_2] \subset {\cal O}_{\varepsilon_1}(A)  \subset {\cal O}_{\varepsilon_1+\eta}(AB)$.
And using remark 0 we achieve
$$[a_1b_1,a_2b_2] \subset {\cal O}_{2\delta} \Bigl([a_1b_1,a_1] \cup [a_1,a_2] \cup [a_2b_2,a_2] \Bigr)
 \subset {\cal O}_{2\delta +max\{ \varepsilon_1+\eta,\varepsilon_3 \}} (AB)~,$$
q.e.d. $\square$

\vspace{.15cm} {\bf \underline{Corollary 1.}} In a hyperbolic
group $G$ every quasiconvex product is a quasiconvex set .

\vspace{.15cm}
This follows directly from the proposition 1 and part (a) of remark 4.

\vspace{.5cm}
{\noindent \bf 3. Intersections of quasiconvex products}

\vspace{.15cm}

Set a partial order on $\Z^2$: $(a_1,b_1) \le (a_2,b_2)$ if $a_1 \le a_2$ and $b_1 \le b_2$. As usual,
$(a_1,b_1) < (a_2,b_2)$ if $(a_1,b_1) \le (a_2,b_2)$ and $(a_1,b_1) \neq (a_2,b_2)$.

\vspace{.15cm}
{\bf Definition :} a finite sequence $\bigl( (i_1,j_1),(i_2,j_2), \dots ,
(i_t,j_t) \bigr)$ of pairs of po\-si\-tive integers
will be called {\it increasing} if it is empty ($t=0$) or
(if $t>0$) $(i_q,j_q) < (i_{q+1}, j_{q+1})$ $\forall ~q=1,2, \dots ,t-1$.
This sequence will also be called $(n,m)$-{\it increasing} ($n,m \in \N$) if
$1 \le i_q \le n$, $1 \le j_q \le m$ for all $q \in \{ 1,2,\dots,t\}$.

\vspace{.15cm}
Note that the length of  an $(n,m)$-increasing sequence never exceeds \\$(n+m-1)$.

Instead of proving theorem 1 we will prove

\vspace{.15cm}
{\bf \underline{Theorem $1'$.}}
Suppose $G_1, \dots, G_n$, $H_1, \dots, H_m$ are quasiconvex sub- \\groups
of the group $G$, $n,m \in \N$; $f,e \in G$.
Then there exist numbers \\ $r, t_l \in \N \cup \{0\}$ and $f_l,\alpha_{lk},
\beta_{lk} \in  G$, $k =1,2,\dots , t_l$ (for every
fixed $l$), \\ $l=1,2,\dots,  r$, such that
$$ fG_1 G_2 \cdot \dots \cdot G_n \cap eH_1 H_2
\cdot \dots \cdot H_m = \bigcup_{l=1}^r f_l S_l \leqno (1)$$
where for each $l$, $t = t_l$, there are indices $1 \le i_1\le i_2\le \dots \le i_t \le n$,
$1 \le j_1\le $ $ \le j_2\le \dots \le j_t \le m$~:
$$ \quad S_l = (G_{i_1}^{\alpha_{l1}} \cap H_{j_1}^{\beta_{l1}}) \cdot
\dots \cdot (G_{i_t}^{\alpha_{lt}} \cap H_{j_t}^{\beta_{lt}}),
\leqno (2) $$
and the sequence $\bigl( (i_1,j_1), \dots, (i_t,j_t)
\bigr)$ is $(n,m)$-increasing.

\vspace{.15cm}
For our convenience, let us also introduce the following

\vspace{.15cm}

{\bf Definition :} the unions as in the right-hand side
of $(1)$ will be called {\it special}. $S_l$ as in $(2)$ will be called
{\it increasing (n,m)-products}.

\vspace{.15cm}
{\bf \underline{Lemma 3.1.}} Consider a geodesic polygon $X_0X_1 \dots X_n$
in the Cayley graph $\Gamma(G,{\cal A})$, $n \ge 2$. Then there are
points $\X_i \in [X_i;X_{i+1}]$, $i=1,2, \dots, n-1$, \\ such that setting
$\X_0 = X_0$, $\X_n = X_n$, we have $(\X_{i-1}|\X_{i+1})_{\X_i} \le \delta$
and $d(\X_i,[\X_{i-1};X_i]) \le \delta$, for $1 \le i \le n-1$.

\vspace{.15cm}
{\it Proof} of the lemma. First, we recursively construct
the vertices $\X_i$. Let \\$\X_1 \in [X_1;X_2],~\U_1 \in [X_0;X_1]$ be
the "special" points of the geodesic triangle $X_0X_1X_2$, i.e.
$|X_1-\X_1| = |X_1-\U_1| = (X_0|X_2)_{X_1}$. Now, if $\X_{i-1}$ is
constructed, denote by $\X_i \in [X_i;X_{i+1}],~\U_i \in [\X_{i-1};X_i]$
the special points of triangle $\X_{i-1}X_iX_{i+1}$
($|X_i-\X_i| = |X_i-\U_i| = (\X_{i-1}|X_{i+1})_{X_i}$). (Figure 2)\\
Then $d(\X_i,[\X_{i-1};X_i]) \le |\X_i-\U_i| \le \delta$, $\forall~
i = 1,2,\dots,n-1$.

\begin{figure}
\begin{center}
\input{fig2.tex}

Figure 2
\end{center}

\end{figure}

For the other part of the claim we will use induction on $n$. \\
$n=2$, then
$$(X_0|X_2)_{\X_1} \stackrel{def}{=} \frac12 (|X_0-\X_1|+|X_2-\X_1|-
|X_0-X_2|) \le $$ $$\le \frac12 (|X_0-\U_1| + |\U_1-\X_1|+|X_2-\X_1|-
|X_0-X_2|) = \frac12 |\U_1-\X_1| \le \frac{\delta}2 \le \delta \quad . $$
Suppose, now, that $n \ge 3$.
Let us evaluate the Gromov product $(\X_0|\X_2)_{\X_1}$.
$$(\X_0|\X_2)_{\X_1} = \frac12 (|X_0-\X_1| + |\X_2-\X_1|-|X_0-\X_2|)~,$$
$|\X_2-\X_1| \le |\X_1-\U_2| + |\X_2-\U_2| \le |\X_1-\U_2| + \delta$,
$|X_0-\X_1| \le |X_0-\U_1| + \delta$, $|X_0-\U_1| + |X_2-\X_1| = |X_0-X_2|$ -- by the definition
of special points of the triangle $X_0X_1X_2$.
Therefore
$$|X_0-\X_1| + |\X_2-\X_1| \le |X_0-\U_1| + |\X_1-\U_2| + 2\delta = $$ $$=
|X_0-\U_1| + (|X_2-\X_1| - |X_2-\U_2|) + 2\delta  = |X_0-X_2| -
|X_2-\X_2| + 2\delta~.$$
Now we notice that $|X_0-X_2|-|X_2-\X_2| \le |X_0-\X_2|$ and obtain:
$$(\X_0|\X_2)_{\X_1} \le \frac12 (|X_0-\X_2| + 2\delta - |X_0-\X_2|) =
\delta~~.$$

To the $n$-gon $\X_1X_2 \dots X_n$ we can apply the induction
hypothesis. \\
The lemma is proved. $\square$

\vspace{.15cm}
{\it Proof} of  theorem $1'$. Define
$T=fG_1 G_2 \cdot \dots \cdot G_n \cap eH_1 H_2
\cdot \dots \cdot H_m$.

Fix some finite generating sets in every $G_i$,$H_j$ and denote
$$K_1 = max\{ 1~;~|f|_G~;~|generators~ of~ G_i|_G~:~i=1,2,\dots,n \} < \infty~,$$
$$K_2 = max\{ 1~;~|e|_G~;~|generators~ of~ H_j|_G~:~j=1,2,\dots,m \} < \infty~.$$

Induction on $(n+m)$.

If $n = 0$ or $m = 0$, then $card(T) \le 1$ and the statement is true.

Let $n \ge 1$ and $ m \ge 1$, $n+m \ge 2$.

\begin{figure}
\begin{center}
\input{fig3.tex}

Figure 3
\end{center}

\end{figure}

Choose an arbitrary $x \in T$, $x = fg_1g_2 \cdot \dots \cdot g_n =
eh_1 h_2 \cdot \dots \cdot h_m$ where $g_i \in G_i$, $h_j \in H_j$,
$i=1,\dots,n$, $j=1,\dots,m$.

Consider a pair of (non-geodesic) polygons associated with $x$ in
$\Gamma(G,{\cal A})$: \\
$P = X_0p_1X_1p_2 \dots p_{n}X_np_0$ and $Q= Y_0q_1Y_1q_2 \dots q_{m}Y_mq_0$
with vertices \\ $X_0 = Y_0=1_{G}$, $X_i=fg_1 \cdot \dots \cdot g_i \in G$, $Y_j = eh_1\cdot \dots \cdot h_j \in G$,
$i=1,\dots ,n$~, $j=1,\dots ,m$, and edges $p_0,p_1,\dots , p_{n}$, $q_0,q_1,\dots,q_{m}$. Such that $p_1$,
starting at $X_0$ and ending at $X_1$,
is a union of a geodesic path corresponding to $f$ and a $G_1$-path corresponding
to $g_1$,
$p_i$ is a $G_i$-path labelled by a word representing the element $g_i$ in $G$ from $X_{i-1}$ to $X_i$,
$i=2,\dots ,n$; $p_0$ is the geodesic path $[X_n,X_0]$ (Figure 3).

By construction, there are constants $\lambda_i,c_i$ (not depending on $x \in T$) such that the segments $p_i$, $i=1,\dots,n$ are
$(\lambda_i,c_i)$-quasigeodesic respectively.

Similarly one constructs the paths $q_j$, $j=0,\dots,m$.

Therefore the geodesic path $p_0=[X_0;X_n] = [Y_0;Y_m]=q_0$ will be labelled
by a word representing $x$ in our Cayley graph.

We will also consider the geodesic polygons $X_0X_1 \dots X_n$ and $Y_0Y_1 \dots Y_m$ with same
vertices as $P$ and $Q$ respectively.


Recalling the property of quasigeodesic paths, for each $i = 1,\dots,n$ \\
 $[~ j = 1,\dots,m~ ]$ we obtain a constant $N_i>0$ [ $M_j>0$ ] (not depending on the element
$x \in T$) such that
$$[X_{i-1},X_{i}] \subset {\cal O}_{N_i} (p_i) \quad  \left[~[Y_{j-1},Y_{j}] \subset {\cal O}_{M_j} (q_j)~ \right] .
\leqno \rm (i)$$

Define $ L= max \{N_1,\dots,N_n,M_1,\dots,M_m\}$.

{\bf a)} Suppose $n,m \ge 2$ (after considering this case, we will see
that the other cases, when $n=1$ or $m=1$ are easier) .

Let's focus our attention on the polygons $X_0 \dots X_n$ and $P$ since
everything for the two others can be done analogously.

One can apply  lemma 3.1 and obtain $\tilde X_i \in
[X_i;X_{i+1}]$, $i=1, \dots,n-1$, such that $(\tilde X_{i-1}|\tilde X_{i+1})_{\tilde X_i} \le
\delta$, $i=1, \dots,n-1$, ($\tilde X_0 = X_0 = 1_G$, $\tilde X_n = X_n = x$), along
with $\tilde U_i \in [\tilde X_{i-1};X_i]$, $|\tilde X_i - \tilde U_i| \le \delta$, $i=1, \dots,n-1$.


\begin{figure}
\begin{center}
\input{fig4.tex}

Figure 4
\end{center}

\end{figure}

Now, using (i), we obtain points $\X_i \in p_{i+1}$, $i=1, \dots,n-1$,
satisfying $d(\tilde X_{i},\X_i) \le L$ ($\X_0 = X_0 = 1_G$, $\X_n = X_n = x$)
and $\U_1 \in p_1$ satisfying $d(\tilde U_1,\U_1) \le L$.
For each $i \in \{ 1,2,\dots,n-2 \}$ the triangle $\tilde X_i \X_i X_{i+1}$
is $\delta$-slim, hence $\exists~ \tilde U'_{i+1} \in [\X_i,X_{i+1}]$~:
$d(\tilde U'_{i+1},\tilde U_{i+1}) \le L+\delta$.
The segment of $p_{i+1}$ between $\X_i$ and $X_{i+1}$ is quasigeodesic with the same constants
as $p_{i+1}$, therefore there is a point $\U_{i+1} \in p_{i+1}$ between $\X_i$ and $X_{i+1}$
such that $d(\tilde U'_{i+1},\U_{i+1}) \le L$, and, consequently, $d(\tilde U_{i+1},\U_{i+1})
\le 2L+\delta$ (see Figure 4).

Let $ \alpha_t$  denote the segment of $p_t$ from $\X_{t-1}$ to $X_{t}$, $t=2,\dots,n$ ,
and \\ $ \beta_s$ -- the subpath of $p_s$ from $\X_{s-1}$ to $\U_{s}$ , $s=1,\dots,n-1$ .
Shifting the points $\X_i, \U_i$ , $i=1,\dots,n-1$ , along their sides of $P$ (~so that $\U_i$ still stays between
$\X_{i-1}$ and $X_i$ on $p_i$~) by
distances at most $K_1$, we can achieve $elem(\beta_1) \in fG_1$
(i.e. $lab(\beta_1)$ represents an element of $fG_1$),
$elem(\alpha_t) \in G_{t+1}$,  $elem(\beta_s) \in G_{s}$, $t=2,\dots,n$, $s=1,2,
\dots,n-1$. And after this,
setting, for brevity, \\ $K = max\{K_1+\frac32 L,K_2+\frac32 L\}$, one obtains
$$(\X_{i-1}|\X_{i+1})_{\X_i} \le \delta + 3 K_1 + 3L \le \delta + 3K \le 14 \delta + 3K
\stackrel{def}{=} C_0~, $$ $$ |\X_i-\U_i| \le \delta + 2K_1+3L+\delta \le 2\delta + 2K ,~i=1,\dots,n-1.$$

Let $elem(\beta_1) = f\g_1$, $elem(\beta_i) = \g_i$, $i=1,\dots,
n-1$, $elem(\alpha_n) = \g_n$, where $\g_k \in G_k$, $k = 1,2,
\dots,n$. $elem([\U_i;\X_i]) = u_i \in G_iG_{i+1}$, $i=1,\dots,n-1$.

Then $|u_i|_G \le 2\delta + 2K$, and there are only finitely many of
possible $u_i$'s for every $i \in \{1,2,\dots,n-1\}$.
Hence, we achieved the following representation for $x$ :
$$ \quad x \stackrel{G}{=} f\g_1 u_1\g_2 u_2\cdot \dots \cdot \g_{n-1} u_{n-1}\g_n~~.
\leqno (*) $$

\begin{figure}
\begin{center}
\input{fig5.tex}

Figure 5
\end{center}

\end{figure}

Similarly, one can obtain
$$ \quad x \stackrel{G}{=} e\h_1 v_1\h_2 v_2\cdot \dots \cdot \h_{m-1} v_{m-1}\h_m~~,
\leqno (**) $$
where $\h_j \in H_j$, $j=1,\dots,m$; $ v_j \in H_jH_{j+1}$ and
$|v_j|_G \le 2\delta + 2K$ for every $j=1,2,\dots,m-1$ (see Figure 5).

${\cal U}_i \stackrel{def}{=} \{ u \in G_iG_{i+1}~:~|u|_G \le 2\delta +2K \}\subset  G$,
$i=1,\dots,n-1~$. $card({\cal U}_i) < \infty$, $\forall ~i = 1,\dots,n-1$.
For convenience, ${\cal U}_0 = {\cal U}_n = G_0 = G_{n+1} \stackrel{def}{=} \{1_G\}$.

Analogously, define ${\cal V}_j \subset H_jH_{j+1}$, $j=1,\dots,m-1$,
and again, \\ ${\cal V}_0 = {\cal V}_m = H_0 = H_{m+1} \stackrel{def}{=} \{1_G\}$.


Set $D=14(\delta + C_0) + 3K = const$, and
${\cal L} = \{ g \in G~:~|g|_G \le D \}$. At last, we denote
$$\Delta_i = {\cal U}_{i-1} \cdot ({\cal L} \cap G_i) \cdot {\cal U}_i
\subset G_{i-1}G_iG_{i+1}\subset G~,~i=1,2,\dots,n~,$$
$$\Theta_i = {\cal V}_{j-1} \cdot ({\cal L} \cap H_j) \cdot {\cal V}_j
\subset H_{j-1}H_jH_{j+1}\subset G~,~j=1,2,\dots,m~.$$

By construction, $card(\Delta_i) < \infty$, $card(\Theta_j) < \infty$,
$\forall ~i,j$.

Take any $i \in \{1,2,\dots,n\}$ and consider the intersection
$$T \supseteq fG_1G_2\cdot \dots \cdot G_{i-1} \Delta_i G_{i+1} \cdot \dots
\cdot G_n \cap eH_1\cdot \dots \cdot H_m = $$
$$ = \bigcup_{g \in \Delta_i} \left[ fG_1G_2\cdot \dots \cdot G_{i-1} g G_{i+1}
\cdot \dots \cdot G_n \cap eH_1\cdot \dots \cdot H_m \right] = $$
$$ = \bigcup_{g \in \Delta_i} \left[ fg(g^{-1}G_1g)(g^{-1}G_2g)\cdot
\dots (g^{-1}G_{i-1}g) G_{i+1} \cdot \dots \cdot G_n \cap eH_1\cdot
\dots \cdot H_m \right] ~.$$
Because of remark 4, one can apply the induction hypothesis to the last expression and conclude
that it is a (finite) "special" union. Hence,
$$T_1 \stackrel{def}{=} \bigcup_{i=1}^n \left(
fG_1G_2\cdot \dots \cdot G_{i-1} \Delta_i G_{i+1}
\cdot \dots \cdot G_n \cap eH_1\cdot \dots \cdot H_m \right) \leqno (3)$$
is also a finite special union.

Because of the symmetry, we parallely showed that
$$T_2 \stackrel{def}{=} \bigcup_{j=1}^m \left( fG_1\cdot \dots \cdot G_n
\cap eH_1H_2\cdot \dots \cdot H_{j-1} \Theta_j H_{j+1} \cdot \dots
\cdot H_m \right) \leqno (4)$$ is a finite "special" union.

We have just proved that there exist $r_1 \in \N \cup \{0\}$, $f_l \in G$ and
increasing $(n,m)$-products $S_l$, $l = 1,2,\dots,r_1$, such that
$$T_1 \cup T_2 = \bigcup_{l=1}^{r_1} f_lS_l \subseteq T~.$$

$T = T_1 \cup T_2 \cup T_3$, where $T_3 \stackrel{def}{=} T\backslash (T_1 \cup T_2)$.
Now, let's consider the case \\ $x \in T_3$. It means that in representations
$(*)$ and $(**)$ for $x$, $|\g_i|_G >D$, \\$|\h_j|_G > D $, for $D= 14(\delta + C_0) + 3K$ and
$\forall~ i=1,\dots,n$, $\forall ~j=1,\dots,m$.\\
Therefore, returning to the pair of polygons we constructed, one
will have~: \\ $||[\X_0;\X_1]||\ge |\g_1|_G - |f|_G - |u_1|_G > 14(\delta +
C_0) + 3K - K - 2\delta - 2K = \\ =12(\delta + C_0) + 2C_0 \stackrel{def}{=} C_1$,
$||[\X_{i-1};\X_i]|| \ge |\g_i|_G - |u_i|_G > 14(\delta + C_0) + 3K - 2\delta -\\
- 2K > C_1$, $i=2,\dots,n-1$,
$||[\X_{n-1};\X_n]|| = |\g_n|_G > C_1$. \\
We also possess the following inequalities~:
$(\X_{i-1}|\X_{i+1})_{\X_i} < C_0$, $i=1, \dots,n-1$, $C_0 \ge 14\delta$,
$C_1 > 12(\delta + C_0)$.

By lemma 1.3, the broken line $[\X_0;\X_1;\dots;\X_n]$ is contained
in the closed $C=2C_0$-neighborhood of the geodesic segment $[\X_0;\X_n]$.
In particular, $$d(\X_{n-1},[\X_0;\X_n]) \le C~. \leqno (5)$$

A similar argument shows that $d(\Y_{m-1},[\Y_0;\Y_m]) \le C$, and, since
$[\X_0;\X_n]= \\ =[X_0;X_n]=[Y_0;Y_m]=[\Y_0;\Y_m]$, one has
$$d(\Y_{m-1},[\X_0;\X_n]) \le C ~. \leqno (6)$$

{\bf b)} In the previous case we assumed that $n,m \ge 2$ and we needed
quite a long argument to prove $(5)$ and $(6)$. On the other hand,
if, for example, $n=1$, then $X_0=\X_{n-1}$ and $(5)$ is trivial.

Because of $(5)$ and $(6)$ one can choose $W,Z \in [X_0;X_n]$ with the
properties $|W-\X_{n-1}| \le C$, $|Z-\Y_{m-1}| \le C$.

The first possibility is, when  the point $W$ on $[X_0;X_n]$ lies between
$Z$ and $X_n$, i.e. $W \in [Z;X_n]$.

\begin{figure}
\begin{center}
\input{fig6.tex}

Figure 6
\end{center}

\end{figure}


Then, since triangles are $\delta$-thin in the hyperbolic space $\Gamma(G,{\cal A})$, \\
$d(W,[\Y_{m-1};X_n]) \le C+ \delta$. Hence $d(\X_{n-1},[\Y_{m-1};X_n]) \le 2C+\delta$. Consequently,
because $q_m$ is quasigeodesic,  there exists a point $R$ on the subpath $\gamma $
of $q_m$ from $\Y_{m-1}$ to $Y_m$
such that  $d(\X_{n-1},R) \le 2C+\delta+K+M_m$ (~$M_m$ is the same as in (i)~)
and $elem([R;Y_m]) = elem(\gamma) ={\hat h}_m \in H_m$.

Define $\Omega = \{g \in G_nH_m~:~|g|_G \le 2C+\delta + K+M_m\}$. Therefore
$card(\Omega) < \infty$ and $elem([\X_{n-1};R]) \in \Omega$.

For each element $g \in \Omega$ take a pair $g'\in G_n$,
$h'\in H_m$ such that $g = g'h'$. By $G'\subset G_n$ denote the set
of all elements $g'$ which we have chosen, by $H'\subset H_m$ -- the set of
all $h'$'s.
$$x = f\g_1u_1 \dots u_{n-1}\g_n = e\h_1v_1 \dots v_{m-1}\h_m~.$$
From the triangle $\X_{n-1}X_nR$ we obtain $\g_n{\hat h}_m^{-1} = g'h' \in
\Omega$, $g' \in G'$, $h' \in H'$. \\
Thus $(g')^{-1}\g_n = h'{\hat h}_m \in G_n \cap H_m$.
$$x \in fG_1G_2 \cdot \dots \cdot G_{n-1} \cdot u_{n-1} g' \cdot
((g')^{-1}\g_n)~\cap~
eH_1H_2 \cdot \dots \cdot H_{m-1} H_m \subset $$

$$\subset fG_1G_2 \cdot \dots \cdot G_{n-1} {\cal U}_{n-1}G'\cdot
(G_n \cap H_m) \cap eH_1H_2 \cdot \dots \cdot H_{m-1} H_m \subset T~.$$
Denote $I = {\cal U}_{n-1}\cdot G'\subset G_{n-1}G_n$ - a finite subset of $G$. Then
$$ x \in fG_1G_2 \cdot \dots \cdot G_{n-1} I\cdot
(G_n \cap H_m) \cap eH_1H_2 \cdot \dots \cdot H_{m-1} H_m = $$
$$=\left[ fG_1G_2 \cdot \dots \cdot G_{n-1} I \cap eH_1H_2 \cdot \dots
\cdot H_m \right] \cdot (G_n \cap H_m) \subset T~.$$

The second possibility, when $Z \in [W;X_n]$ is considered analogously,
and, in this case, one obtains a finite subset $J \subset H_{m-1}H_m$ such that
$$ x \in \left[ fG_1G_2 \cdot \dots \cdot G_n \cap eH_1H_2 \cdot \dots
\cdot H_{m-1} J \right] \cdot (G_n \cap H_m) \subset T~.$$

Therefore, we showed that \quad
$T_3 \subseteq \left[ T'_3 \cup T''_3 \right] \cdot (G_n \cap H_m) \subset T$ \quad where
$$T'_3 \stackrel{def}{=} fG_1G_2 \cdot \dots \cdot G_{n-1} I \cap eH_1H_2 \cdot \dots
\cdot H_m  ~, \leqno (7)$$
$$T''_3 \stackrel{def}{=}  fG_1G_2 \cdot \dots \cdot G_n \cap eH_1H_2 \cdot \dots
\cdot H_{m-1} J ~. \leqno (8)$$

Combining the formulas $(3)$,$(4)$,$(7)$,$(8)$ and the property that if
$H \le G$ and $a \in H$ then $aH=Ha=H$, we obtain the following

\vspace{.15cm} \underline{\bf Lemma 3.2.} In  notations of the
theorem 1
$$fG_1G_2 \cdot \dots \cdot G_{n} \cap eH_1H_2 \cdot \dots \cdot H_m = T_1 \cup T_2 \cup
 \left[ T'_3 \cup T''_3 \right] \cdot (G_n \cap H_m) $$
where
$$T_1 = \bigcup_{i=1}^n \left(
fG_1G_2\cdot \dots \cdot G_{i-1} \bar \Delta_i G_{i+1}
\cdot \dots \cdot G_n \cap eH_1\cdot \dots \cdot H_m \right)~,$$
$$T_2 = \bigcup_{j=1}^m \left( fG_1\cdot \dots \cdot G_n
\cap eH_1H_2\cdot \dots \cdot H_{j-1} \bar \Theta_j H_{j+1} \cdot \dots
\cdot H_m \right)~,$$
$$T'_3 = fG_1G_2 \cdot \dots \cdot G_{n-1} \bar I \cap eH_1H_2 \cdot \dots
\cdot H_m  ~,$$
$$T''_3 =  fG_1G_2 \cdot \dots \cdot G_n \cap eH_1H_2 \cdot \dots
\cdot H_{m-1} \bar J ~ $$
for some finite subsets $\bar \Delta_i \subset G_i$,$\bar \Theta_j \subset H_j$,
$\bar I \subset G_n$,$\bar J \subset H_m$, $1\le i \le n$, $1 \le j \le m$.

\vspace{.15cm}
Now, to finish the proof of the theorem, we apply the inductive hypothesis:
$$T_3 \subseteq \bigcup_{g \in I}\left[ fG_1G_2 \cdot \dots \cdot G_{n-1} g
\cap eH_1H_2 \cdot \dots \cdot H_m \right] \cdot (G_n \cap H_m) ~\cup$$
$$~\cup \bigcup_{h \in J}\left[ fG_1G_2 \cdot \dots \cdot G_n \cap eH_1H_2
\cdot \dots \cdot H_{m-1} h \right] \cdot (G_n \cap H_m) = $$
$$= \bigcup_{g \in I}\left[ fgG_1^{g^{-1}}G_2^{g^{-1}} \cdot \dots \cdot
G_{n-1}^{g^{-1}} \cap eH_1H_2 \cdot \dots \cdot H_m \right] \cdot
(G_n \cap H_m)~ \cup$$
$$~\cup \bigcup_{h \in J}\left[ fG_1G_2 \cdot \dots \cdot G_n \cap eh
H_1^{h^{-1}} H_2^{h^{-1}} \cdot \dots \cdot H_{m-1}^{h^{-1}} \right]
\cdot (G_n \cap H_m) = $$
$$= \left( \bigcup_{g \in I}\left[ \bigcup_{k=1}^{\tilde r} {\tilde f}_k
{\tilde S}_k \right] \cup \bigcup_{h \in J}\left[ \bigcup_{q=1}^{\hat r}
{\hat f}_q {\hat S}_q \right] \right) \cdot (G_n \cap H_m) = $$
$$= \bigcup_{l=r_1+1}^r f_lS_l \subset T~.$$
Here ${\tilde r}, {\hat r}, r \in \N \cup \{0\}$, $r \ge r_1$,
${\tilde f}_k,{\hat f}_q,f_l \in G$; ${\tilde S}_k$ is an (n-1,m)-increasing
product, ${\hat S}_q$ is an (n,m-1)-increasing product and
$S_l$ is an (n,m)-increasing product; \\ $k=1,\dots,{\tilde r}$; $q=1\dots,
{\hat r}$; $l = r_1+1,\dots,r$.

Hence,
$$ T = T_1 \cup T_2 \cup T_3 \subseteq \bigcup_{l=1}^r f_lS_l \subseteq T~,$$
and, thus
$$T = \bigcup_{l=1}^r f_lS_l~.$$

So, the theorem is proved. $\square$


\vspace{.15cm}
{\it Proof} of Corollary 2. Observe that arbitrary quasiconvex product $f_1G_1f_2G_2 \cdot \dots \cdot f_nG_n$
is equal to a "transformed" product $f G'_1 G'_2 \cdot \dots \cdot G'_n$ where \\
$G'_i = (f_{i+1} \cdot \dots \cdot f_n)^{-1}G_i(f_{i+1} \cdot \dots \cdot f_n)$, $i=1,\dots,n-1$,
$G'_n = G_n$, are quasiconvex subgroups of $G$ by remark 4 and $f = f_1f_2\cdot \dots f_n \in G$.
It remains to apply theorem 1 to the intersection of "transformed products" several times
because a $(n,m)$-increasing product is also a quasiconvex product. $\square$

\vspace{.5cm}
{ \bf \noindent 4. Products of elementary subgroups }

\vspace{.15cm}
Recall that a group $H$ is called {\it elementary} if it has a
cyclic subgroup $\langle h \rangle$ of finite index.

\vspace{.15cm}
{\bf \underline{Remark 5.}}   An elementary subgroup $H$ of a hyperbolic group $G$
is quasiconvex .

\vspace{.15cm}
Indeed, we have : $|H:\langle h \rangle | <\infty $ .
If the element $h$ has a finite order , then $H$ is finite and, thus, quasiconvex. In the case,
when the order of $h$ is infinite, by lemmas 1.2,1.1 $\langle h \rangle$ is a quasiconvex subgroup of
$G$. By remark 4 and lemma 2.1 $H$ is quasiconvex.

It is well known that any element $x$ of infinite order in $G$ is contained in a
unique maximal elementary subgroup $E(x) \leqslant G$ (see [4]). And the intersection
of two distinct maximal elementary subgroups in a hyperbolic group is finite. Any infinite
elementary subgroup contains an element of infinite order.

Obviously, a conjugate subgroup to a maximal elementary subgroup is also maximal
elementary.



\vspace{.15cm}
{\it Proof} of Theorem 2. The sufficiency is trivial.

Without loss of generality one can assume $n \ge m$. In this case theorem 2 immediately follows from

\vspace{.15cm}
{\bf \underline{Theorem $2'$}} Let $n \ge m$, $G_1,G_2,\dots,G_n,H_1,H_2,\dots,H_m$ be infinite maximal elementary
subgroups of $G$, $f,e \in G$, and $g_i \in G_i$, $i=1,2,\dots,n$, be elements of infinite order.
Also, assume $G_i \neq G_{i+1}$, $H_j \neq H_{j+1}$, $i=1,\dots,n-1$, $j=1,\dots,m-1$. If there is a
sequence of positive integers $(t_k)_{k=1}^\infty$ with the properties:
$$\lim_{k\to \infty} t_k = \infty \mbox{~~and~~}
fg_1^{t_k}g_2^{t_k} \cdot \dots \cdot g_n^{t_k} \in eH_1H_2\cdot \dots \cdot H_m~ \mbox{~for all~} k\in \N,$$ then
$n=m$, $G_n=H_n$, and there exist elements $z_i \in H_i$, $i=1,\dots,n$, such that
$G_i = (z_{n} z_{n-1} \dots z_{i+1}) \cdot H_i \cdot
(z_{n} z_{n-1} \dots z_{i+1})^{-1}$, $i=1,2,\dots,n-1$, $f = ez_1^{-1}z_2^{-1}\dots z_n^{-1}$.
Consequently, $fG_1\cdot \dots \cdot G_n=eH_1 \cdot \dots \cdot H_m$.

\vspace{.15cm}
In the conditions of theorem $2'$, let $h_j \in H_j$ be fixed elements of
infinite order, $j=1,2,\dots,m$. Then  $G_i=E(g_i)$, $H_j=E(h_j)$ and $|G_i:\langle g_i \rangle | < \infty$,
$|H_j:\langle h_j \rangle | <  \infty$.
Hence, there exists $T \in \N$ such that for all $j$ and $\forall~v \in H_j$ ~~$\exists~ \beta \in \Z$, $y \in H_j$~:
$v = y \cdot h^{\beta}_j$ and $|y|_G \le T$. Thus, every element $h \in eH_1\cdot \dots \cdot H_m$
can be presented in the form
$$ h= ey_1h^{\beta_1}_1 y_2h^{\beta_2}_2 \cdot \dots \cdot y_m h^{\beta_m}_m           \leqno (9)$$
where $\beta_j \in \Z$, $y_j \in H_j$, $|y_j|_G \le T$, $j=1,2,\dots, m$.

\vspace{.15cm}
{\bf Definition:} the representation (9) for $h$ will be called {\it reduced} if for any $i,j$, $1\le i< j \le m$,
such that $\beta_i, \beta_j \neq 0$, one has
$$(y_{i+1}h^{\beta_{i+1}}_{i+1}\dots h^{\beta_{j-1}}_{j-1}y_j)^{-1}\cdot h_i \cdot (y_{i+1}h^{\beta_{i+1}}_{i+1}\dots h^{\beta_{j-1}}_{j-1}y_j)
\notin H_j = E(h_j)~.$$

\vspace{.15cm}
Observe that each element $h \in eH_1\cdot \dots \cdot H_m$ has a reduced representation. Indeed, if
$(y_{i+1}h^{\beta_{i+1}}_{i+1}\dots h^{\beta_{j-1}}_{j-1}y_j)^{-1}\cdot h_i \cdot (y_{i+1}h^{\beta_{i+1}}_{i+1}\dots h^{\beta_{j-1}}_{j-1}y_j)
\in H_j $ for some $1\le i< j \le m$ then there are $\beta'_j \in \Z$, $y'_j \in H_j$, $|y'_j|_G\le T$~:
$$y_j \cdot (y_{i+1}h^{\beta_{i+1}}_{i+1}\dots h^{\beta_{j-1}}_{j-1}y_j)^{-1}\cdot h^{\beta_i}_i \cdot (y_{i+1}h^{\beta_{i+1}}_{i+1}\dots h^{\beta_{j-1}}_{j-1}y_j)
\cdot h^{\beta_j}_j = y'_jh^{\beta'_j}_j ~.$$
Therefore,
$$h= ey_1h^{\beta_1}_1 \cdot \dots \cdot y_{i-1} h^{\beta_{i-1}}_{i-1} y_i y_{i+1} h^{\beta_{i+1}}_{i+1} \cdot \dots \cdot
y_{j-1}h^{\beta_{j-1}}_{j-1} y'_jh^{\beta'_j}_j y_{j+1}h^{\beta_{j+1}}_{j+1} \cdot \dots \cdot y_m h^{\beta_m}_m $$
and the number of non-zero $\beta_k$'s is decreased. Continuing this process, we will obtain a reduced representation for $h$
after a finite number of steps .

\vspace{.15cm}
{\it Proof} of Theorem $2'$. Let $h_j \in H_j$, $1 \le j \le m$, $T$, be as above.
Induction on $n$.

If n=1, then, evidently, $m=1$, and $\forall~ k \in \N$ there is $y_{t_k} \in H_1$,
$|y_{t_k}|_G \le T$, and $\beta_{t_k} \in \Z$ such that $fg_1^{t_k} = ey_{t_k}h_1^{\beta_{t_k}}$.
Because of having $\lim_{k\to \infty} t_k = \infty$, one can choose $p,q \in \N$ so that $t_p < t_q$ and
$y_{t_p} = y_{t_q}$. Therefore,
$$fg_1^{t_p}h_1^{-\beta_{t_p}} =  ey_{t_p} = fg_1^{t_q}h_1^{-\beta_{t_q}}~, $$
and, thus, $g_1^{t_p-t_q} = h_1^{\beta_{t_p}-\beta_{t_q}}$ -- an element of infinite order in the intersection
of $G_1$ and $H_1$. Consequently, $G_1=H_1$, because these subgroups are maximal elementary.

Assume, now, that $n>1$.
For every $k \in \N$ one has
$$ fg^{t_k}_1g^{t_k}_2 \cdot \dots \cdot g^{t_k}_n = ey_{t_k1}h^{\beta_{t_k1}}_1 y_{t_k2}h^{\beta_{t_k2}}_2
\cdot \dots \cdot y_{t_km} h^{\beta_{t_km}}_m \leqno (10)$$
where the product in the right-hand side is reduced. Obviously,
there exists a subsequence $(l_k)_{k=1}^\infty $ of $(t_k)$ and $C \in \N$ such
that for each $j \in \{ 1,2,\dots,m\}$ either  $|\beta_{l_k j}| \le C$ for all $k$
or $\lim_{k \to \infty} |\beta_{l_k j}| = \infty$. \\
Therefore, since $|y_{l_k j}|_G \le T$ $\forall~k \in \N$, $\forall~j $,
there is a subsequence $(s_k)_{k=1}^\infty$ of $(l_k)$ such that $y_{s_k j} = y_j \in H_j$ $\forall~j$, and
if  for $ j \in \{1,\dots,m\}$ we had $|\beta_{l_k j}| \le C$ $\forall~k \in \N$ then
$|\beta_{s_k j}| = \beta_j \in \Z$ $\forall~k\in \N$,
and $\lim_{k \to \infty} |\beta_{s_k j}| = \infty$ for all other $j$'s.

Thus, $\{1,2,\dots,m\} = J_1 \cup J_2$ where if $j \in J_1$ then $|\beta_{s_k j}|= \beta_j$ for every $k$, and if $j \in J_2$
then $\lim_{k \to \infty} |\beta_{s_k, j}| = \infty$. Let $J_2 = \{j_1,j_2, \dots, j_{\varkappa}\} \subset \{1,2,\dots,m\}$,
$j_1 < j_2 < \dots <j_\varkappa$,
and denote
$$w_1 = y_{1}^{-1} \in H_1 \mbox{ if } j_1 = 1~, \mbox{otherwise, if $j_1 > 1$, } $$
$$w_1 =y_{j_1}^{-1}h_{j_1-1}^{-\beta_{j_1-1}} y_{j_1-1}^{-1}\cdot \dots \cdot
h_{1}^{-\beta_{1}}y_{1}^{-1} \in H_{j_1}H_{j_1-1}\cdot \dots \cdot H_1~;$$
$$\dots \dots $$
$$w_\varkappa = y_{j_{\varkappa}}^{-1} \in H_{j_\varkappa} \mbox{ if } j_\varkappa = j_{\varkappa-1} +1~,
\mbox{otherwise, if $  j_\varkappa > j_{\varkappa-1} +1$,}  $$
$$w_{\varkappa} =y_{j_{\varkappa}}^{-1}h_{j_{\varkappa}-1}^{-\beta_{j_{\varkappa}-1}} y_{j_{\varkappa}-1}^{-1}
\cdot \dots \cdot h_{j_{\varkappa-1}+1}^{-\beta_{j_{\varkappa-1}+1}}y_{j_{\varkappa-1}+1}^{-1}
\in H_{j_\varkappa}H_{j_{\varkappa-1}} \cdot \dots \cdot H_{j_{\varkappa-1}+1}~;$$
$$w_{\varkappa+1} = 1_G \mbox{ if } j_{\varkappa}=m ~, \mbox{ otherwise, if }
j_{\varkappa} < m, $$ $$w_{\varkappa+1} =h_{m}^{-\beta_{m}} y_{m}^{-1}\cdot \dots \cdot
h_{j_{\varkappa}+1}^{-\beta_{j_{\varkappa}+1}}y_{j_{\varkappa}+1}^{-1} \in H_m H_{m-1} \cdot \dots
\cdot H_{j_{\varkappa}+1}~.$$

To simplify the formulas, denote $\delta_{k \nu} = -\beta_{s_k,j_{\nu}}$, $1\le \nu \le \varkappa$.\\
Then $\lim_{k \to \infty} |\delta_{k \nu}|= \infty$ for every $\nu=1,2,\dots,\varkappa$.
(10) is equivalent to
$$ u_k \stackrel{def}{=} f g_1^{s_k} g^{s_k}_2 \cdot \dots \cdot g^{s_k}_{n-1} w_{\varkappa+1}
h^{\delta_{k\varkappa}}_{j_{\varkappa}} w_{\varkappa} h^{\delta_{k,\varkappa-1}}_{j_{\varkappa-1}}
\cdot \dots \cdot w_{2} h^{\delta_{k1}}_{j_1}w_1 e^{-1} =1_G \leqno (11) $$

So, $|u_k|_G = 0$ for all $k \in \N$. Denote $K = max\{|f|_G,|w_{1}e^{-1}|_G,|w_{2}|_G,$
$\dots,|w_{\varkappa+1}|_G \}$,    and assume that
$g_n \notin w_{\varkappa+1}E(h_{j_{\varkappa}}) w^{-1}_{\varkappa+1}$. The product in the right-hand side
of (10) was reduced, therefore
$h_{j_{\nu}} \notin w_\nu E(h_{j_{\nu-1}}) w_{\nu-1}^{-1}$, $\nu = 2,3,\dots,\varkappa$.
Thus, we can apply Lemma 1.8 to (11) and obtain $\lambda >0$, $c \ge 0$ and $M>0$
(depending on $K$, $g_1,\dots,g_n$,$h_{j_{1}},\dots,h_{j_\varkappa}$) such that
if $s_k \ge M$ and $|\delta_{k\nu}| \ge M$, $\nu=2,3,\dots,\varkappa$, then $|u_k|_G \ge \lambda \cdot s_k -c$.
Now, by the choice of the sequence
$(s_k)$, there exists $N \in \N$~: $s_k > M$ and $|\delta_{k\nu}| > M$ $\forall~k \ge N$, $\nu=2,3,\dots,\varkappa$.
Thus, taking $k \ge max\{N,c/\lambda\}+1$, we achieve a contradiction: $0 = |u_k|_G < \lambda \cdot s_k -c$.

Hence, $g_n \in w_{\varkappa+1}E(h_{j_{\varkappa}}) w^{-1}_{\varkappa+1}$ which implies
$$G_n=E(g_n)=w_{\varkappa+1}E(h_{j_{\varkappa}}) w^{-1}_{\varkappa+1}=
E(w_{\varkappa+1}h_{j_{\varkappa}} w^{-1}_{\varkappa+1}). \leqno (12)$$ Consequently, for every $k \in \N$ \quad
$w_{\varkappa+1}^{-1}g_n^{s_k} w_{\varkappa+1}h_{j_{\varkappa}}^{\delta_{k,\varkappa}} = y_{k j_{\varkappa}}'
h_{j_{\varkappa}}^{\gamma_k} \in H_{j_\varkappa}$
where $|y_{k j_{\varkappa}}'|_G \le T$. By passing to a subsequence of $(s_k)$ we can assume that
\\ $y_{k j_{\varkappa}}' = y_{j_\varkappa}' \in H_{j_\varkappa}$ for every $k$.
Therefore
$$ u_k = f g_1^{s_k} g^{s_k}_2 \cdot \dots \cdot g^{s_k}_{n-1} w_{\varkappa+1} y'_{j_\varkappa}
h^{\gamma_k}_{j_\varkappa} w_{\varkappa} h^{\delta_{k,\varkappa-1}}_{j_{\varkappa-1}}
\cdot \dots \cdot w_{2} h^{\delta_{k1}}_{j_1}w_1 e^{-1} =1_G~. $$

Suppose $\limsup_{k \to \infty} |\gamma_k| = \infty$. Since $E(g_{n-1}) = G_{n-1}  \neq G_n = E(g_n)$, we have
$g_{n-1} \notin w_{\varkappa+1}E(h_{j_{\varkappa}}) w^{-1}_{\varkappa+1} =
w_{\varkappa+1} y_{k j_{\varkappa}}' E(h_{j_{\varkappa}}) (y_{k j_{\varkappa}}')^{-1} w^{-1}_{\varkappa+1}$
(~because $y_{k j_{\varkappa}}' \in H_{j_\varkappa} =  E(h_{j_{\varkappa}})$~).

Then for $K' = max\{K,|w_{\varkappa+1}y_{j_{\varkappa}}'|_G\}$ by Lemma 1.8 there exist
$\lambda >0$, $c \ge 0$ and $M>0$
(depending on $K'$, $g_1,\dots,g_{n-1}$,$h_{j_{1}},\dots,h_{j_\varkappa}$) such that
if $s_k \ge M$, $|\delta_{k\nu}| \ge M$, $\nu=2,3,\dots,\varkappa$, and $|\gamma_k| \ge M$ then
$|u_k|_G \ge \lambda \cdot s_k -c$. Now, by the assumption on $(s_k)$ and $(\gamma_k)$,
there exists $N \in \N$, $N>c/\lambda$, such that $s_N > M$, $|\delta_{N\nu}| > M$,
$\nu=2,3,\dots,\varkappa$, and $|\gamma_N|>M$. Which leads us to a contradiction: $0 = |u_k|_G < \lambda \cdot s_k -c$.

Thus, $|\gamma_k| \le C_1$ for some constant $C_1$, so,
by passing to a subsequence as above, we can assume that $\gamma_k = \gamma$
$\forall~k \in \N$. Hence, after setting \\ $z_\varkappa = w_{\varkappa+1}y'_{j_\varkappa} h^{\gamma}_{j_\varkappa}
w_{\varkappa}$, for every natural index $k$ we will have
$$ u_k = f g_1^{s_k} g^{s_k}_2 \cdot \dots \cdot g^{s_k}_{n-1} z_\varkappa
h^{\delta_{k,\varkappa-1}}_{j_{\varkappa-1}} w_{j_{\varkappa-1}}
\cdot \dots \cdot w_{2} h^{\delta_{k1}}_{j_1}w_1 e^{-1} =1_G ~. $$
Which implies $f g_1^{s_k} g^{s_k}_2 \cdot \dots \cdot g^{s_k}_{n-1} \in
ew_1^{-1} H_{j_1} w_2^{-1} H_{j_2} \cdot \dots \cdot w_{j_{\varkappa-1}}^{-1} H_{j_{\varkappa-1}} z_{\varkappa}^{-1}
= $ $=u H_1^{v_2} H_2^{v_3} \cdot \dots \cdot H_{j_{\varkappa-1}}^{v_{\varkappa}}$ where
$v_\nu = z_{\varkappa}w_{\varkappa-1} \cdot \dots \cdot w_{\nu+1} w_{\nu}$, $\nu=2,3,\dots,\varkappa-1$,
$v_{\varkappa} = z_{\varkappa}$, $u = e w_1^{-1} w_2^{-1} \cdot \dots \cdot w_{\varkappa-1}^{-1}z_{\varkappa}^{-1}$.

$n-1 \ge m-1 \ge \varkappa-1 $ and the other conditions of the theorem $2'$ are satisfied, therefore
one can apply the induction hypothesis and obtain that $n-1 = \varkappa-1$, hence, $\varkappa=m=n$,
$j_\nu = \nu$, $1\le \nu \le \varkappa$, and, by definition, $w_\nu=y_{j_\nu}^{-1} \in H_\nu$, $\nu=1,2,\dots,n$,
$w_{\varkappa+1} = 1_G$, $z_\varkappa = z_n \in H_n$.  And also
$G_{n-1} = H_{j_{\varkappa-1}}^{v_{\varkappa}} = H_{n-1}^{z_n}$, and there exist ${\hat z}_i \in H_i$,
$1\le i \le n-1$, such that
$$G_i = ({\hat z}_{n-1}^{v_n} {\hat z}_{n-2}^{v_{n-1}} \dots
{\hat z}_{i+1}^{v_{i+2}}) \cdot H_i^{v_{i+1}} \cdot
({\hat z}_{n-1}^{v_n} {\hat z}_{n-2}^{v_{n-1}} \dots
{\hat z}_{i+1}^{v_{i+2}})^{-1} = $$
$$ = (z_nz_{n-1} \cdot \dots \cdot z_{i+1}) \cdot H_i \cdot (z_nz_{n-1} \cdot \dots \cdot z_{i+1})^{-1}
~,~i=1,2,\dots,n-2~,$$ where $z_p = {\hat z}_p w_p \in H_p$, $1 \le p \le n-1$,
$f = u \left({\hat z}_1^{v_2} \right)^{-1} \left({\hat z}_2^{v_3} \right)^{-1} \cdot \dots \cdot
\left({\hat z}_{n-1}^{v_2} \right)^{-1} =$ $= ez_1^{-1}z_2^{-1} \cdot \dots \cdot z_n^{-1}$.

By (12) $G_n = E(h_n) = H_n$.
The proof of the theorem $2'$ is finished. $\square$

\vspace{.15cm}
Suppose $G_1,G_2,\dots,G_n$ are {\it infinite maximal elementary} subgroups of $G$,
$f_1,\dots,f_n \in G$, $n \in \N \cup \{0 \}$.

\vspace{.15cm}
{\bf Definition :} the set $P = f_1G_1f_2G_2\cdot \dots \cdot f_nG_n$ will be called {\it ME-product}.
Thus, if $n=0$, we have the empty set. For convenience, we will also consider every element $g \in G$  to be
a ME-product.
As in the proof of corollary 2, every such ME-product can be brought to a form (however, not unique)
$$P' = fG_1'G_2'\cdot \dots G_k'$$ where $0 \le k \le n$,
$f \in G$, $G_i'$ are infinite maximal elementary subgroups, $i=1,2,\dots,k$, and $G_i' \neq G_{i+1}'$,
$1 \le i\le k-1$. The number $k$ in this case will be called {\it rank} of the ME-product $P$
(thus, $rank(P)=rank(P')=k \le n$). \\
A set $U$ which can be presented as a finite
union of ME-products has rank $k$, by definition, if $U = \bigcup_{i=1}^t P_i$~, where
$P_i$~, $i=1,\dots,t$, are ME-products, and $ k=max\{ rank(P_i)~|~1\le i \le t \}$~.

\vspace{.15cm}
{\it Note}: an empty set is defined to have rank $(-1)$; any element of the group $G$ is a
ME-product of rank $0$; thus
any finite non-empty subset of $G$ is a finite union of ME-products of rank $0$.

\vspace{.15cm}
{\bf \underline{Remark 6.}}  the rank of a ME-product is defined correctly by theorem 2.
By theorem $2'$ the definition of the rank of a finite union of ME-products is correct.

\vspace{.15cm}
{\bf \underline{Lemma 4.1.}} Suppose  $P$,$R$ are ME-products in a hyperbolic group  $G$. Then the intersection
$T \stackrel{def}{=} P \cap R $ is a finite union of ME-products and its rank is at most $rank(P)$.
If $rank(T)=rank(P)$ then $T=P$~.

\vspace{.15cm}
{\it Proof. } Since a conjugate to an infinite maximal elementary subgroup is also infinite maximal elementary,
it follows from theorem 1 that $T$ is a finite union of ME-products $P_i$, $1\le i \le t$
(for some $t \in \N \cup \{0\}$):
  $$T=P \cap R = \bigcup_{i=1}^t P_i~.$$

For each $i=1,\dots,t$, $P_i \subseteq P$, therefore by theorem $2'$, $rank(P_i) \le rank(P)$
(otherwise we would get a contradiction), and $rank(P_i) = rank(P)$ if and only if $P_i=P$.
Thus $ rank(T)=max\{ rank(P_i)~|~1\le i \le t \} \le rank(P)$. If $rank(T) = rank(P)$ then $rank(P_i)=rank(P)$ for
some $i$, and so, $P_i=P=T$. \quad Q.e.d. $\square$

\vspace{.15cm}
\noindent As an immediate consequence of lemma 4.1 one obtains

\vspace{.15cm}
\underline{\bf Corollary 3.} let $P$ be a ME-product of rank $n$ and  $U$ be a finite union
of ME-products. Then the set $P\cap U$ is a finite union of ME-products, \\
$rank(P \cap U) \le n$~, and   if $rank(P \cap U) = n$ then $P \cap U = P$~.

\vspace{.15cm}
\underline{\bf Corollary 4.} A non-elementary hyperbolic group $G$ can not be equal to a finite union
of its ME-products.

\vspace{.15cm}
{\it Proof.} Suppose, by the contrary, that $G$ is a finite union of ME-products:
$G=P_1 \cup \dots \cup P_l$ and $rank(G)=m$. Since $G$ is
not elementary, there exist two elements $x,y \in G$ of infinite order such that $E(x) \neq E(y)$. Hence, one
can construct a ME-product $P=G_1G_2 \cdot \dots \cdot G_{m+1}$ in $G$ where $G_i=E(x)$ if $i$ is even, and
$G_i=E(y)$ if $i$ is odd. Consequently, $rank(P) = m+1$, but $P \subset G$, thus
$$P \cap G = P= \bigcup_{j=1}^l \left( P_j \cap P \right)~.$$
By lemma 4.1, $rank(P_j \cap P) \le rank(P_j) \le m$ ~~for every $j = 1,2,\dots, l$. Therefore, we achieve a contradiction
with the definition of rank : $m+1 = $ $=rank(P) = rank(P \cap G) \le m$. $\square$

\vspace{.15cm}
A group $H$ is called {\it bounded-generated} if it is a product of finitely many cyclic subgroups, i.e.
there are elements $x_1,x_2,\dots,x_k \in H$ such that every $h \in H$ is equal to
$x_1^{s_1}x_2^{s_2}\cdot \dots \cdot x_k^{s_k}$ for some $s_1,\dots,s_k \in \Z$.

\vspace{.15cm}
\underline{\bf Corollary 5.} Any bounded-generated hyperbolic group is elementary.

\vspace{.15cm}
{\it Proof.} Indeed, any cyclic subgroup of a hyperbolic group either is finite or
is contained in some infinite maximal elementary subgroup. Hence, their product is contained in a finite
union of ME-products and we can apply corollary 4. $\square$

\vspace{.15cm}
{\it Proof} of Theorem 3. Since there exist at most countably many different ME-products in $G$, it is enough to
consider  only their countable intersections.
Let $P_{ji}$, $1 \le i \le k_j$, $k_j,j \in \N$, be ME-products, and
$U_j = \bigcup_{i=1}^{k_j} P_{ji}$~-- their finite unions. Let
$$T = \bigcap_{j=1}^\infty U_j~.$$
One has to show that there exist ME-products $R_1,\dots,R_s$, $s\in \N \cup \{0\}$,
such that \quad $T= R_1 \cup \dots \cup R_s$~.

Induct on $n = rank(U_1)$.
$$T =  \left( \bigcup_{i=1}^{k_1} P_{1i} \right) \cap \bigcap_{j=2}^\infty U_j =
\bigcup_{i=1}^{k_1} \left( P_{1i}  \cap \bigcap_{j=2}^\infty U_j \right)
$$
So, it is enough to consider the case when $k_1 = 1$~, $U_1 = P_{11} = P$. \\
If $n = 0$ then $P$ is finite and there is nothing to prove. \\
Assume that $n>0$ and let $J \in \N$ be the smallest index such that \\
$P\cap U_{J} \neq P$~(if there is no such $J$ then $T = P$ and the
theorem is true). Therefore

$$T = P \cap \bigcap_{j=J}^\infty U_j = \left( P \cap U_J \right) \cap \bigcap_{j=J+1}^\infty U_j ~.$$

By corollary 3, $P \cap U_J$ is a finite union of ME-products~:
$$P \cap U_J  = \bigcup_{l=1}^t R'_l~, ~t \in \N \cup \{0\}$$
 and $rank(P \cap U_J)<n$ because of the choice of $J$, therefore $rank(R'_l)<n$~, \\
$\forall~l=1,2,\dots,t$.

Hence, by the induction hypothesis,
$$T = \bigcup_{l=1}^t \left[ R'_l \cap \bigcap_{j=J+1}^\infty U_j \right] =  \bigcup_{l=1}^t \left[
R_{l1} \cup \dots \cup R_{ls_l} \right]$$
for some ME-products $R_{l1},\dots,R_{ls_l}$, $s_l\in \N\cup \{0\}$, $1 \le l \le t$~.
$\square$

\vspace{.15cm}
The statement of the theorem 3 fails to be true if
maximal elementary subgroups in the definition of ME-products
one substitutes by arbitrary elementary subgroups.
Below we construct an example to demonstrate that .

Let $G = F(x,y)$ be the free group with two generators,
$q_1<q_2<q_3< \dots$ ~~be an infinite sequence of prime numbers . Define
$d_i = q_1 q_2 \cdot \dots \cdot q_i$, \\$c_i = q_1q_2 \cdot \dots \cdot q_{i-1}q^2_i = d_i \cdot q_i$, $i \in \N$, and
the sets
$P_i$, $i \in \N$, as follows~:  \\ $P_1= \langle x^{d_1}  \rangle $ -- cyclic subgroup of $G$
generated by $x^{d_1}=x^{q_1}$~, \\
$P_2 = \langle y \rangle \cdot \langle y x^{c_1} y^{-1} \rangle \cdot \langle y^2 x^{d_2} y^{-2}
\rangle\cdot \langle y \rangle $~, \\
$P_3 = \langle y \rangle \cdot \langle y x^{c_1} y^{-1} \rangle \cdot \langle y^2 x^{c_2} y^{-2}
\rangle\cdot  \langle y^3 x^{d_3} y^{-3} \rangle\cdot \langle y \rangle$~, \\
$\dots \dots$ \\
$P_i = \langle y \rangle \cdot \langle y x^{c_1} y^{-1} \rangle \cdot \langle y^2 x^{c_2} y^{-2}
\rangle\cdot \dots \cdot \langle y^{i-1} x^{c_{i-1}} y^{-(i-1)} \rangle
 \langle y^i x^{d_i} y^{-i}\rangle\cdot \langle y \rangle$~,\\
$\dots \dots$

Now consider the intersection \quad $T = \bigcap_{i=1}^\infty P_i $~.
Let us observe that \\ $P_1 \cap P_2 = \langle x^{c_1} \rangle \cup  \langle x^{d_2} \rangle$, $\dots$,
$\bigcap_{i=1}^k P_i =  \langle x^{c_1} \rangle \cup \dots \cup \langle x^{c_{k-1}} \rangle \cup \langle x^{d_k} \rangle$, $\dots$~.
\\Indeed, $P_1 \cap P_2 = \langle x^{d_1} \rangle  \cap \left( \langle x^{c_1} \rangle \cup \langle x^{d_2} \rangle \right) = $
$\langle x^{c_1} \rangle \cup  \langle x^{d_2} \rangle$. Inducting on $k$, we get
$$\bigcap_{i=1}^k P_i = \left( \bigcap_{i=1}^{k-1} P_i \right) \cap P_k =
\langle x^{c_1} \rangle \cup \dots \cup \langle x^{c_{k-2}} \rangle \cup \langle x^{d_{k-1}} \rangle \cap
\left( \langle x^{c_1} \rangle \cup \dots \cup \langle x^{c_{k-1}} \rangle \right. \cup $$
$$\cup \left. \langle x^{d_{k}} \rangle \right)   =\langle x^{c_1} \rangle \cup \dots \cup \langle x^{c_{k-2}} \rangle \cup
\langle x^{d_{k-1}} \rangle \cap \left( \langle x^{c_{k-1}} \rangle \cup \langle x^{d_{k}} \rangle \right) = $$
$$=\langle x^{c_1} \rangle \cup \dots \cup \langle x^{c_{k-1}} \rangle \cup \langle x^{d_{k}} \rangle~.$$
Since $\bigcap_{i=1}^\infty \langle x^{d_{i}} \rangle= \{1\}$, therefore $T = \bigcup_{i=1}^\infty \langle x^{c_i} \rangle$~.

If $q_1=2,q_2=3,q_3=5,\dots$, is chosen to be the enumeration of all primes, one can show directly that the set $T$ can
not be presented as a finite union of products $f_1G_1f_2G_2\cdot \dots \cdot f_nG_n$~, where $f_1, \dots, f_n \in G$ and
$G_1, \dots, G_n$ are elementary (in this case cyclic) subgroups of $G$~.
We are not going to
do that, instead we will use a set-theoretical argument~: there are only countably many such
finite unions, hence there is an infinite sequence of primes $q_1<q_2<q_3< \dots$  such that
the corresponding set $\bigcap_{i=1}^\infty P_{i} $  is the example sought (~because
the sets $\bigcap_{i=1}^\infty P_{i}$ and $\bigcap_{i=1}^\infty P'_{i}$ corresponding to different
increasing sequences of prime numbers $\alpha = \{ q_1,q_2,q_3, \dots \}$ and  $\alpha'=\{q'_1,q'_2,q'_3, \dots\}$
are distinct:
if $q_l \in \alpha \backslash \alpha'$ then $x^{c_l} \in \bigcup_{i=1}^\infty \langle x^{c_i} \rangle \backslash
\bigcup_{j=1}^\infty \langle x^{c'_j} \rangle$~)~.

\vspace{1.5cm}
{\Large \bf \noindent Acknowledgements}

The author is grateful to his advisor Professor
A.Yu. Ol'shanskii  for suggesting problems, helpful discussions and comments.

\end{document}